\documentclass{amsart}
\usepackage{graphicx}
\usepackage{amscd}
\usepackage{amsmath}
\usepackage{amsfonts}
\usepackage{amssymb}
\newtheorem{theorem}{Theorem}
\theoremstyle{plain}

\newtheorem{corollary}{Corollary}

\newtheorem{lemma}{Lemma}

\numberwithin{equation}{section}

\begin{document}
\author{}
\thanks{M\'{e}canique et Gravitation, Universit\'{e} de Mons-Hainaut, Facult\'{e} des
Sciences, 15 Avenue Maistriau, B-7000 Mons, Belgium}
\thanks{Universit\`{a} degli Studi di Bergamo, Facolt\`{a} di Ingegneria, \noindent
Viale Marconi, 5, 24044 Dalmine \noindent(Bergamo) Italy }
\thanks{E-mail: Garattini@mi.infn.it }
\title{}
\maketitle

\begin{center}
\bigskip{\large {\textbf{REMO GARATTINI}}} \vspace{1cm}

\textrm{\textbf{{\large {Harnack's Inequality on Homogeneous Spaces}}}}\vspace{1cm}
\end{center}

\noindent\textbf{Abstract.} We consider a homogeneous space $X=\left(
X,d,m\right)  $ of dimension $\nu\geq1$ and a local regular Dirichlet form in
$L^{2}\left(  X,m\right)  .$ We prove that if a Poincar\'{e} inequality holds
on every pseudo-ball $B\left(  x,R\right)  $ of $X$, then an Harnack's
inequality can be proved on the same ball with local characteristic constant
$c_{0}$ and $c_{1}$

\section{Introduction and Results}

We consider a \textit{connected, locally compact topological space }$X$. We
suppose that a \textit{distance} $d$ is defined on $X$ and we suppose that the
balls
\[
B\left(  x,r\right)  =\left\{  y\in X:d\left(  x,y\right)  <r\right\}  ,\qquad
r>0,
\]
form a basis of open neighborhoods of $x\in X.$ Moreover, we suppose that a
(positive) Radon measure $m$ is given on $X$, with supp$m=X.$ The triple
$\left(  X,d,m\right)  $ is assumed to satisfy the following property: There
exist some constants $0<R_{0}\leq+\infty,\nu>0$ and $c_{0}>0,$ such that
\begin{equation}
0<c_{0}\left(  \frac{r}{R}\right)  ^{\nu}m\left(  B\left(  x,R\right)
\right)  \leq m\left(  B\left(  x,r\right)  \right)  \label{a1}%
\end{equation}
for every $x\in X$ and every $0<r\leq R<R_{0}.$ Such a triple $\left(
X,d,m\right)  $ will be called a homogeneous space of dimension $\nu.$ We
point out, however, that a given exponent $\nu$ occurring in $\left(
\ref{a1}\right)  $ should be considered, more precisely, as an upper bound of
the ``homogeneous dimension'', hence we should better call $\left(
X,d,m\right)  $ a homogeneous space of dimension less or equal than $\nu.$ We
further suppose that we are also given a \textit{strongly local, regular,
Dirichlet form a }in the Hilbert space $L^{2}\left(  X,m\right)  $ - in the
sense of M. Fukushima \cite{fuku}, - whose domain in $L^{2}\left(  X,m\right)
$ we shall denote by $\mathcal{D}\left[  a\right]  $. Furthermore, we shall
restrict our study to Dirichlet forms of diffusion type, that is to forms
\textit{a} that have the following \textit{strong local property:} $a\left(
u,v\right)  =0$ for every $u,v\in\mathcal{D}\left[  a\right]  $ with $v$
constant on supp $u.$ We recall that the following integral representation of
the form \textit{a }holds
\[
a\left(  u,v\right)  =\int_{X}\mu\left(  u,v\right)  \left(  dx\right)
,\qquad u,v\in\mathcal{D}\left[  a\right]  ,
\]
where $\mu\left(  u,v\right)  $ is a uniquely defined signed Radon measure on
X, such that $\mu\left(  d,d\right)  \leq m$, with $d\in\mathcal{D}%
_{loc}\left[  a\right]  $: this last condition is fundamental for the
existence of cut off functions associated to the distance. Moreover, the
restriction of the measure $\mu\left(  u,v\right)  $ to any open subset
$\Omega$ of $X$ depends only on the restrictions of the functions $u,v$ to
$\Omega.$ Therefore, the definition of the measure $\mu\left(  u,v\right)  $
can be unambiguously extended to all $m$-measurable functions $u,v$ on $X$
that coincide $m-a.e.$ on every compact subset of $\Omega$ with some functions
of $\mathcal{D}\left[  a\right]  .$ The space of all such functions will be
denoted by $\mathcal{D}_{loc}\left[  a,\Omega\right]  .$ Moreover we denote by
$\mathcal{D}\left[  a,\Omega\right]  $ the closure of $\mathcal{D}\left[
a\right]  \cap C_{0}\left[  \Omega\right]  $ in $\mathcal{D}\left[  a\right]
$. The homogeneous metric $d$ and the energy form $a$ associated to the
\textit{energy measure }$\mu$, both given on a relatively compact open subset
$X_{0}$ of $X$ with $\overline{\Omega}\subset X_{0}$, are then assumed to be
mutually related by the following basic assumption:

There exists an exponent $s>2$ and constants $c_{1}$, $c_{1}>0$ and $k\geq1$,
such that \cite{Biroli}:

\begin{description}
\item [ i)]for every $x\in X_{0}$ and every $0<R<R_{0}$ the following
Poincar\'{e} inequality holds:
\begin{equation}
\int\limits_{B\left(  x,R\right)  }\left|  u-\bar{u}_{B\left(  x,R\right)
}\right|  ^{2}dm\leq c_{1}R^{2}\int\limits_{B\left(  x,kR\right)  }\mu\left(
u,u\right)  \left(  dx\right)  \label{a2}%
\end{equation}
for all $u\in D\left[  a,B\left(  x,kR\right)  \right]  $, where
\[
\bar{u}_{B\left(  x,R\right)  }=\frac{1}{m\left(  B\left(  x,R\right)
\right)  }\int\limits_{B\left(  x,R\right)  }udm.
\]

\item[ ii)] By i), for every $x\in X_{0}$ and every $0<R<R_{0}$, the following
Sobolev inequality of exponent $s$ can be proved:
\begin{equation}
\left(  \frac{1}{m\left(  B\left(  x,R\right)  \right)  }\int\limits_{B\left(
x,R\right)  }\left|  u\right|  ^{s}dm\right)  ^{\frac{1}{s}}\leq c_{1}R\left(
\int\limits_{B\left(  x,R\right)  }\mu\left(  u,u\right)  \left(  dx\right)
\right)  ^{\frac{1}{2}},\label{a3}%
\end{equation}
where $u\in\mathcal{D}\left[  a,B\left(  x,kR\right)  \right]  $ and supp
$u\subset B\left(  x,R\right)  $. Instead of working with the fixed constants
of $\left(  \ref{a1}\right)  $ and$\left(  \ref{a2},\ref{a3}\right)  $, we
make this simple generalization
\begin{equation}
c_{0}\rightarrow c_{0}\left(  x\right)  \qquad\text{and\qquad}c_{1}\rightarrow
c_{1}\left(  r\right)  ,\label{aa1}%
\end{equation}
where $c_{0}\left(  x\right)  ,c_{0}^{-1}\left(  x\right)  \in L_{loc}%
^{\infty}\left(  X_{0}\right)  $ and $c_{1}\left(  r\right)  $ is a decreasing
function of $r$. Then we assume:

\item[ iii)]
\begin{equation}
\int\limits_{B\left(  x,R\right)  }\left|  u-\bar{u}_{B\left(  x,R\right)
}\right|  ^{2}dm\leq c_{1}^{2}\left(  R\right)  R^{2}\int\limits_{B\left(
x,kR\right)  }\mu\left(  u,u\right)  \left(  dx\right)  .\label{a1a}%
\end{equation}
By iii) the following Sobolev type inequality may be proved

\item[ iv)]
\begin{equation}
\left(  \frac{1}{m\left(  B\left(  x,R\right)  \right)  }\int\limits_{B\left(
x,R\right)  }\left|  u\right|  ^{s}dm\right)  ^{\frac{1}{s}}\leq\tau^{3}%
c_{1}\left(  R\right)  R\left(  \frac{1}{m\left(  B\left(  x,R\right)
\right)  }\int\limits_{B\left(  x,R\right)  }\mu\left(  u,u\right)  \left(
dx\right)  \right)  ^{\frac{1}{2}},\label{a2a}%
\end{equation}
\end{description}

where $u$ is as in i) and ii) and where we have defined $\tau=\left(
\underset{B\left(  x,2R\right)  }{\sup}\frac{1}{c_{0}\left(  x\right)
}\right)  ^{\frac{1}{2}}$. Our purpose will be the Harnack's inequality
recovery for Dirichlet forms when the substitution $\left(  \ref{aa1}\right)
$ is performed.

\begin{theorem}
[Harnack]\label{t1}Let $\left(  \ref{a1}\right)  $ , $\left(  \ref{a1a}%
\right)  $, $\left(  \ref{a2a}\right)  $ and $\left(  \ref{aa1}\right)  $
hold, and let $u$ be a non-negative solution of $\ a\left(  u,v\right)  =0$.
Let $\mathcal{O}$ be an open subset of $X_{0}$ and$\ u\in\mathcal{D}%
_{loc}\left[  \mathcal{O}\right]  ,\ \forall v\in\mathcal{D}_{0}\left[
\mathcal{O}\right]  $ with $B\left(  x,r\right)  \subset\mathcal{O}$, then
\[
\underset{B_{\frac{1}{2}}}{ess\ \sup}u\leq\exp\gamma\mu^{\prime}\mu
^{2}\ \underset{B_{\frac{1}{2}}}{ess\ \inf\ u},
\]
where $\gamma$ $\equiv\gamma\left(  \nu,k\right)  $, $k$ a positive constant,
$\mu^{\prime}=\tau c_{1}\left(  \frac{r}{2}\right)  $and $\mu=\tau^{3}%
c_{1}\left(  \frac{r}{2}\right)  $. A standard consequence of the previous
Theorem is the following
\end{theorem}

\begin{corollary}
Suppose that
\begin{equation}
\int\limits_{r}^{R}e^{-\gamma\mu\left(  x,\rho\right)  }\frac{d\rho}{\rho
}\rightarrow\infty\text{\qquad for\qquad}r\rightarrow0\label{aaa1}%
\end{equation}
then the solution is continuous in the point under consideration. In
particular, if $\mu\left(  x,\rho\right)  \approx o\left(  \log\log\frac
{1}{\rho}\right)  $, then there exists $c>0$ such that
\[
\underset{B\left(  x,r\right)  }{osc}u\leq c\frac{\left(  \log\frac{1}%
{R}\right)  }{\left(  \log\frac{1}{r}\right)  }\underset{B\left(  x,R\right)
}{osc}u.
\]
\end{corollary}

From the point of view of partial differential equations these results can be
applied to two important classes of operators on $\mathbb{R}^{n}$:

\begin{description}
\item [ a)]\textit{Doubly Weighted uniformly elliptic operators} in divergence
form with measurable coefficients, whose coefficient matrix $A=\left(
a_{ij}\right)  $ satisfies
\[
w\left(  x\right)  \left|  \xi\right|  ^{2}\leq\left\langle A\xi
,\xi\right\rangle \leq v\left(  x\right)  \left|  \xi\right|  ^{2}.
\]
Here $\left\langle \cdot,\cdot\right\rangle $ denotes the usual dot product;
$w$ and $v$ are weight functions, respectively belonging to $A_{2}$ and
$D_{\infty}$ such that the following \textit{Poincar\'{e}} inequality
\[
\left(  \frac{1}{\left|  v\left(  B\right)  \right|  }\int\limits_{B}\left|
f\left(  x\right)  -f_{B}\right|  ^{q}vdx\right)  ^{\frac{1}{2}}\leq cr\left(
\frac{1}{\left|  w\left(  B\right)  \right|  }\int\limits_{B}\left|  \nabla
f\right|  ^{2}wdx\right)  ^{\frac{1}{2}}%
\]
holds.
\end{description}

\begin{description}
\item \textit{Doubly Weighted H\"{o}rmander type operators}\cite{franchi}$,$
whose form is $L=X_{k}^{\ast}\left(  \alpha^{hk}\left(  x\right)
X_{h}\right)  $ where $X_{h},h=1,\ldots,m$ are smooth vector fields in
$\mathbb{R}^{n}$ that satisfy the H\"{o}rmander condition and $\alpha=\left(
\alpha^{hk}\right)  $ is any symmetric $m\times m$ matrix of measurable
functions on $\mathbb{R}^{n},$such that
\[
w\left(  x\right)  \sum\limits_{i}\left\langle X_{i},\xi\right\rangle ^{2}%
\leq\sum\limits_{i,j}\alpha_{ij}\left(  x\right)  \xi_{i}\xi_{j}\leq v\left(
x\right)  \sum\limits_{i}\left\langle X_{i},\xi\right\rangle ^{2},
\]
where $X_{i}\xi\left(  x\right)  =\left\langle X_{i},\nabla\xi\right\rangle
,\ i=1,\ldots,m$, $\nabla\xi$ is the usual gradient of $\xi$ and $\left\langle
\cdot,\cdot\right\rangle $ denotes the usual inner product on $\mathbb{R}%
^{n}.$ Then the following \textit{Poincar\'{e}} inequality for vector fields
\[
\left(  \frac{1}{\left|  v\left(  B\right)  \right|  }\int\limits_{B}\left|
f\left(  x\right)  -f_{B}\right|  ^{q}vdx\right)  ^{\frac{1}{2}}\leq cr\left(
\frac{1}{\left|  w\left(  B\right)  \right|  }\int\limits_{B}\left(
\sum\limits_{j}\left|  \left\langle X_{j},\nabla f\left(  x\right)
\right\rangle \right|  ^{2}\right)  ^{\frac{1}{2}}wdx\right)  ^{\frac{1}{2}},
\]
holds, with $w\in A_{2}$ and $v$ $\in D_{\infty}$.
\end{description}

\section{Harnack's Inequality}

We prove the Harnack's Inequality of Theorem $\left(  \ref{t1}\right)  $ for a
non negative solution $u\geq\delta>0$ and with a constant $C$ independent of
$\delta.$ The result of Theorem $\left(  \ref{t1}\right)  $ is obtained
passing to the limit $\delta\rightarrow0.$

\begin{lemma}
\label{l1}Assume that $\left(  \ref{a1}\right)  ,$ $\left(  \ref{a1a}\right)
$, $\left(  \ref{a2a}\right)  $ and $\left(  \ref{aa1}\right)  $ hold.

Let $u$ be a non-negative subsolution of $a\left(  u,v\right)  =0,\ u\in
\mathcal{D}_{loc}\left[  \mathcal{O}\right]  ,\ \forall v\in\mathcal{D}%
_{0}\left[  \mathcal{O}\right]  $. Define $\tau=\left(  \underset{B\left(
x,2R\right)  }{\sup}\frac{1}{c_{0}\left(  x\right)  }\right)  ^{\frac{1}{2}}$
then
\[
\left(  \underset{B_{\alpha}}{ess\sup}u\right)  ^{p}\leq\left(  c\tau^{3}%
\frac{c_{1}\left(  \frac{r}{2}\right)  }{\left(  1-\alpha\right)  }\right)
^{2\frac{\sigma}{\sigma-1}}\left(  \frac{1}{m\left(  B\right)  }%
\int\limits_{B}u^{p}m\left(  dx\right)  \right)  ,
\]
where $\sigma=\frac{q}{2}$, $p\geq2$ and
\[
q=\left\{
\begin{array}
[c]{cc}%
\frac{2\nu}{\nu-2} & \text{if }\nu>2\\
\text{any value} & \text{if }\nu\leq2
\end{array}
\right.
\]
\end{lemma}%

\proof
We prove the result for a bounded non negative subsolutions. Let $\beta\geq1$
and $0<M<\infty$, define $H_{M}\left(  t\right)  =t^{\beta}$ for $t\in\left[
0,M\right]  $ and $H_{M}\left(  t\right)  =M^{\beta}+\beta M^{\beta-1}\left(
t-M\right)  $ for $t>M.$ For fixed $M$ define
\[
\phi_{k}\left(  x\right)  =\eta\left(  x\right)  ^{2}\int\limits_{0}%
^{u_{k}\left(  x\right)  }H_{M}^{^{\prime}}\left(  t\right)  ^{2}dt
\]
and
\[
a\left(  u,\phi\right)  :=\int\limits_{X}\mu\left(  u,\phi\right)  \left(
dx\right)  =\int\limits_{X}\mu\left(  u,\eta\left(  x\right)  ^{2}%
\int\limits_{0}^{u}H_{M}^{^{\prime}}\left(  t\right)  ^{2}dt\right)  \left(
dx\right)  \leq0,
\]
then
\begin{equation}%
\begin{array}
[c]{c}%
\int\limits_{X}\mu\left(  u,\eta\left(  x\right)  ^{2}\int\limits_{0}^{u}%
H_{M}^{^{\prime}}\left(  t\right)  ^{2}dt\right)  \left(  dx\right)
=\int\limits_{X}\mu\left(  u,\eta\left(  x\right)  ^{2}\right)  \int
\limits_{0}^{u}H_{M}^{^{\prime}}\left(  t\right)  ^{2}dt\left(  dx\right)  \\
\\
+\int\limits_{X}\mu\left(  u,\int\limits_{0}^{u}H_{M}^{^{\prime}}\left(
t\right)  ^{2}dt\right)  \eta\left(  x\right)  ^{2}\left(  dx\right)  =\\
\\
\int\limits_{X}2\eta\mu\left(  u,\eta\right)  \int\limits_{0}^{u}%
H_{M}^{^{\prime}}\left(  t\right)  ^{2}dt\left(  dx\right)  +\int
\limits_{X}\mu\left(  u,u\right)  H_{M}^{^{\prime}}\left(  u\right)  ^{2}%
\eta\left(  x\right)  ^{2}\left(  dx\right)  \leq0,
\end{array}
\label{ll1}%
\end{equation}
and therefore
\begin{equation}
\int\limits_{X}\mu\left(  u,u\right)  H_{M}^{^{\prime}}\left(  u\right)
^{2}\eta\left(  x\right)  ^{2}\left(  dx\right)  \leq2\int\limits_{X}\left|
\eta\mu\left(  u,\eta\right)  \int\limits_{0}^{u}H_{M}^{^{\prime}}\left(
t\right)  ^{2}dt\right|  \left(  dx\right)  .\label{ll2}%
\end{equation}
Taking account of the inequality
\[
2\left|  fg\right|  \left|  \mu\left(  u,v\right)  \right|  \leq f^{2}%
\mu\left(  u,u\right)  +g^{2}\mu\left(  v,v\right)
\]
we get
\[
2\left|  \mu\left(  u,\eta\right)  \right|  \left|  \eta\int\limits_{0}%
^{u}H_{M}^{^{\prime}}\left(  t\right)  ^{2}dt\right|  \left(  dx\right)
\leq\frac{1}{2}\mu\left(  u,u\right)  \eta_{{}}^{2}H_{M}^{^{\prime}}\left(
u\right)  ^{2}+2\mu\left(  \eta,\eta\right)  \left[  \frac{1}{H_{M}^{^{\prime
}}\left(  u\right)  }\int\limits_{0}^{u}H_{M}^{^{\prime}}\left(  t\right)
^{2}dt\right]  ^{2}%
\]
and putting into$\left(  \ref{ll1}\right)  $ we have
\[
\frac{1}{2}\int\limits_{X}\mu\left(  u,u\right)  H_{M}^{^{\prime}}\left(
u\right)  ^{2}\eta\left(  x\right)  ^{2}\left(  dx\right)  \leq2\int
\limits_{X}\mu\left(  \eta,\eta\right)  \left[  \frac{1}{H_{M}^{^{\prime}%
}\left(  u\right)  }\int\limits_{0}^{u}H_{M}^{^{\prime}}\left(  t\right)
^{2}dt\right]  ^{2}\left(  dx\right)
\]%
\begin{equation}
\leq2\int\limits_{X}\mu\left(  \eta,\eta\right)  \left[  uH_{M}^{^{\prime}%
}\left(  u\right)  \right]  ^{2}\left(  dx\right)  .\label{ll3}%
\end{equation}
Let us consider $\eta$ as a cut-off function s.t. for $\frac{1}{2}\leq s<t<1$
$B\left(  x,sr\right)  \subset B\left(  x,tr\right)  \subset\vartheta
\subset\subset X,\eta\equiv0$ on $X-B\left(  x,tr\right)  ,\eta\equiv1$ on
$B\left(  x,sr\right)  ,0\leq\eta<1$ on $X.$ Then
\[%
\begin{array}
[c]{c}%
\int\limits_{B_{s}}\mu\left(  H_{M}\left(  u\right)  ,H_{M}\left(  u\right)
\right)  \eta\left(  x\right)  ^{2}\left(  dx\right)  \leq\int
\limits_{\vartheta}\mu\left(  H_{M}\left(  u\right)  ,H_{M}\left(  u\right)
\right)  \eta\left(  x\right)  ^{2}\left(  dx\right)  \\
\\
\leq\frac{40}{\left(  t-s\right)  ^{2}r^{2}}\int\limits_{B_{t}}\left[
uH_{M}^{^{\prime}}\left(  u\right)  \right]  ^{2}m\left(  dx\right)  .
\end{array}
\]
With the result of $\left(  \ref{a2}\right)  $, applied to $H_{M}\left(
u\right)  $, one gets
\[%
\begin{array}
[c]{c}%
\left(  \frac{1}{m\left(  B_{s}\right)  }\int\limits_{B_{s}}\left|
H_{M}\left(  u\right)  -\bar{H}_{M}\left(  u\right)  _{B_{s}}\right|
^{q}m\left(  dx\right)  \right)  ^{\frac{1}{q}}\leq c\tau^{2}c_{1}\left(
sr\right)  sr\left[  \frac{1}{m\left(  B_{s}\right)  }\int\limits_{B_{s}}%
\mu\left(  H_{M}\left(  u\right)  ,H_{M}\left(  u\right)  \right)  \right]
^{\frac{1}{2}}\\
\\
\leq c\tau^{2}c_{1}\left(  sr\right)  sr\left(  \frac{40}{\left(  t-s\right)
^{2}r^{2}}\frac{1}{m\left(  B_{s}\right)  }\int\limits_{B_{t}}\left[
uH_{M}^{^{\prime}}\left(  u\right)  \right]  ^{2}m\left(  dx\right)  \right)
^{\frac{1}{2}}\\
\\
\leq c\tau^{2}\sqrt{40}c_{1}\left(  sr\right)  \frac{s}{\left(  t-s\right)
}\left(  \frac{m\left(  B_{t}\right)  }{m\left(  B_{s}\right)  }\right)
^{\frac{1}{2}}\left(  \frac{1}{m\left(  B_{t}\right)  }\int\limits_{B_{t}%
}\left[  uH_{M}^{^{\prime}}\left(  u\right)  \right]  ^{2}m\left(  dx\right)
\right)  ^{\frac{1}{2}}\\
\\
\leq c\tau^{3}\sqrt{40}c_{1}\left(  sr\right)  \frac{s}{\left(  t-s\right)
}\left(  \frac{t}{s}\right)  ^{\frac{\nu}{2}}\left(  \frac{1}{m\left(
B_{t}\right)  }\int\limits_{B_{t}}\left[  uH_{M}^{^{\prime}}\left(  u\right)
\right]  ^{2}m\left(  dx\right)  \right)  ^{\frac{1}{2}}.
\end{array}
\]
The average of $H_{M}\left(  u\right)  $ on $B_{s}$ is defined by:
\[%
\begin{array}
[c]{c}%
\underset{B_{s}}{av}H_{M}\left(  u\right)  =\frac{1}{m\left(  B_{s}\right)
}\int\limits_{Bs}H_{M}\left(  u\right)  m\left(  dx\right)  \leq\tau
^{2}\left(  \frac{t}{s}\right)  ^{\nu}\frac{2}{m\left(  B_{t}\right)  }%
\int\limits_{B_{s}}H_{M}\left(  u\right)  m\left(  dx\right)  \\
\\
\leq\tau\left(  \frac{t}{s}\right)  ^{\frac{\nu}{2}}\left(  \frac{1}{m\left(
B_{t}\right)  }\int\limits_{B_{t}}\left[  uH_{M}^{^{\prime}}\left(  u\right)
\right]  ^{2}m\left(  dx\right)  \right)  ^{\frac{1}{2}}%
\end{array}
\]
Recall that $m\left(  B_{t}\right)  \leq\tau^{2}m\left(  B_{s}\right)  \left(
\frac{t}{s}\right)  ^{\nu}$ and $H_{M}\left(  u\right)  \leq uH_{M}^{^{\prime
}}\left(  u\right)  .$ Therefore
\[%
\begin{array}
[c]{c}%
\left(  \frac{1}{m\left(  B_{s}\right)  }\int\limits_{B_{s}}\left|
H_{M}\left(  u\right)  \right|  ^{q}dm\right)  ^{\frac{1}{q}}=\left(  \frac
{1}{m\left(  B_{s}\right)  }\int\limits_{B_{s}}\left|  H_{M}\left(  u\right)
-\bar{H}_{M}\left(  u\right)  _{B_{s}}+\bar{H}_{M}\left(  u\right)  _{B_{s}%
}\right|  ^{q}dm\right)  ^{\frac{1}{q}}\\
\\
\leq c\left[  \left(  \frac{1}{m\left(  B_{s}\right)  }\int\limits_{B_{s}%
}\left|  H_{M}\left(  u\right)  -\bar{H}_{M}\left(  u\right)  _{B_{s}}\right|
^{q}dm\right)  ^{\frac{1}{q}}+\left(  \frac{1}{m\left(  B_{s}\right)  }%
\int\limits_{B_{s}}\left|  \bar{H}_{M}\left(  u\right)  _{B_{s}}\right|
^{q}dm\right)  ^{\frac{1}{q}}\right]  \\
\\
\leq c\tau^{3}\left(  \frac{t}{s}\right)  ^{\frac{\nu}{2}}\left(  \frac
{c_{1}\left(  sr\right)  s}{\left(  t-s\right)  }+1\right)  \left(  \frac
{1}{m\left(  B_{t}\right)  }\int\limits_{B_{t}}\left[  uH_{M}^{^{\prime}%
}\left(  u\right)  \right]  ^{2}m\left(  dx\right)  \right)  ^{\frac{1}{2}};
\end{array}
\]
where in the last expression we have included all irrelevant constants into
$c$. From the fact that $\frac{s}{\left(  t-s\right)  }+1\leq2\frac{s}{\left(
t-s\right)  },$ $\frac{s}{\left(  t-s\right)  }\geq1.$ Therefore due to
$uH_{M}^{^{\prime}}\left(  u\right)  \leq u\beta u^{\beta-1}=\beta u^{\beta
},H_{M}\left(  u\right)  \geq u_{\chi_{\left\{  \tilde{u}\leq M\right\}  }%
}^{\beta},$ we obtain by letting $M\rightarrow\infty$ that
\[
\left(  \frac{1}{m\left(  B_{s}\right)  }\int\limits_{B_{s}}u^{\beta
q}dm\right)  ^{\frac{1}{q}}\leq c\tau^{3}\left(  \frac{t}{s}\right)
^{\frac{\nu}{2}}\left(  \frac{\beta c_{1}\left(  sr\right)  s}{\left(
t-s\right)  }\right)  \left(  \frac{1}{m\left(  B_{t}\right)  }\int
\limits_{B_{t}}u^{2\beta}m\left(  dx\right)  \right)  ^{\frac{1}{2}}.
\]
Raise both sides to the power $\frac{1}{\beta}$ and putting $2\beta=\tilde{r}$
and $q=2\sigma$ we have
\[
\left(  \frac{1}{m\left(  B_{s}\right)  }\int\limits_{B_{s}}u^{\sigma
r}dm\right)  ^{\frac{1}{\sigma\tilde{r}}}\leq c\left[  \tau^{3}\left(
\frac{t}{s}\right)  ^{\frac{\nu}{2}}\left(  \frac{\beta c_{1}\left(
sr\right)  s}{\left(  t-s\right)  }\right)  \right]  ^{\frac{2}{\tilde{r}}%
}\left(  \frac{1}{m\left(  B_{t}\right)  }\int\limits_{B_{t}}u^{\tilde{r}%
}m\left(  dx\right)  \right)  ^{\frac{1}{\tilde{r}}}.
\]
Now, starting from fixed $\alpha$ and $p,\frac{1}{2}\leq\alpha<1,p\geq2$
iterate this inequality for $t$ and $s$ successive entries in the sequence
$s_{j}=\alpha+\frac{\left(  1-\alpha\right)  }{\left(  j+1\right)
},j=0,1\ldots,$ and $\tilde{r}$ and $\sigma\tilde{r}$ successive entries in
$\left\{  \sigma^{j}p\right\}  ,$ recalling that $\sigma>1.$ If we denote
$a_{j}=\frac{s_{j}}{s_{j+1}-s_{j}}$ the conclusion of the lemma is a
consequence of the estimate of
\[
\log\prod_{j=0}^{\infty}\left[  c\tau^{3}c_{1}\left(  s_{j}r\right)  \left(
\frac{s_{j+1}}{s_{j}}\right)  ^{\frac{\nu}{2}}\left(  \sigma^{j}pa_{j}\right)
\right]  ^{\frac{2}{\sigma^{j}p}}=\sum_{j=0}^{\infty}\frac{2}{\sigma^{j}p}%
\log\left[  c\tau^{3}c_{1}\left(  s_{j}r\right)  \left(  \sigma^{j}%
pa_{j}\right)  \right]
\]
where we have used the fact that $\frac{s_{j+1}}{s_{j}}=\frac{\alpha
+\frac{\left(  1-\alpha\right)  }{\left(  j+2\right)  }}{\alpha+\frac{\left(
1-\alpha\right)  }{\left(  j+1\right)  }}\simeq1.$ Therefore
\[%
\begin{array}
[c]{c}%
\sum_{j=0}^{\infty}\frac{2}{\sigma^{j}p}\log\left(  c\tau^{3}c_{1}\left(
s_{j}r\right)  p\right)  +\sum_{j=0}^{\infty}\frac{2}{\sigma^{j}p}\log\left(
\sigma^{j}a_{j}\right)  \\
\\
\leq\frac{2}{p}\frac{\sigma}{\sigma-1}\log\left(  c\tau^{3}c_{1}\left(
\frac{r}{2}\right)  p\right)  +\sum_{j=0}^{\infty}\frac{2}{\sigma^{j}p}%
\log\sigma^{j}+\sum_{j=0}^{\infty}\frac{2}{\sigma^{j}p}\left[  2\log\left(
j+2\right)  +\log\frac{1}{1-\alpha}\right]  \\
\\
\leq\frac{2}{p}\frac{\sigma}{\sigma-1}\log\left(  c\tau^{3}c_{1}\left(
\frac{r}{2}\right)  p\right)  +\sum_{j=0}^{\infty}\frac{2}{\sigma^{j}p}%
\log\sigma^{j}+\frac{2}{p}\frac{\sigma}{\sigma-1}\log\frac{1}{1-\alpha}%
+4\sum_{j=0}^{\infty}\frac{\log\left(  j+2\right)  }{\sigma^{j}p}\\
\\
\leq\frac{2}{p}\frac{\sigma}{\sigma-1}\log\left(  \frac{cp}{1-\alpha}%
c_{1}\left(  \frac{r}{2}\right)  \tau^{3}\right)
\end{array}
\]
and at the end
\[
\prod_{j=0}^{\infty}\left[  c\tau^{3}c_{1}\left(  s_{j}r\right)  \left(
\frac{s_{j+1}}{s_{j}}\right)  ^{\frac{\nu}{2}}\left(  \sigma^{j}pa_{j}\right)
\right]  ^{\frac{2}{\sigma^{j}p}}\leq\left(  cc_{1}\left(  \frac{r}{2}\right)
\frac{\tau^{3}}{1-\alpha}\right)  ^{\frac{2}{p}\frac{\sigma}{\sigma-1}},
\]
where we used the following inequalities $\log\left(  a_{j}\right)  \leq
2\log\left(  j+2\right)  +\log\frac{1}{1-\alpha},$ $p^{\frac{1}{p}}\leq c$ for
$p\geq2.$ The general case can be obtained if we apply the above result for
$p=2$ to the truncated subsolutions and we obtain that $u$ is locally bounded.
Coming back to the previous result, we obtain the Lemma $\left(
\ref{l1}\right)  $ by an approximation by truncation.

\begin{lemma}
\label{l2} With the same notation and hypothesis as in Lemma $\left(
\ref{l1}\right)  $ and $-\infty<p<+\infty,$
\end{lemma}%

\[
\left(  \underset{B_{\alpha}}{ess\sup}u^{p}\right)  \leq\left(  \frac
{c\tau^{3}c_{1}\left(  \frac{r}{2}\right)  p+1}{\left(  1-\alpha\right)
}\right)  ^{\frac{2\sigma}{\left(  \sigma-1\right)  }}\left(  \frac
{1}{m\left(  B_{t}\right)  }\int_{B_{t}}u^{p}m\left(  dx\right)  \right)  .
\]%
\proof
It is sufficient to consider $-\infty<p<2.$ Define $\phi\left(  x\right)
=\eta^{2}\left(  x\right)  u^{\beta}\left(  x\right)  ,$ $\eta\left(
x\right)  \geq0,$ $-\infty<\beta<+\infty.$ Then
\[%
\begin{array}
[c]{c}%
0\geq\int_{X}\mu\left(  u,\phi\right)  \left(  dx\right)  =\int_{X}\mu\left(
u,\eta^{2}u^{\beta}\right)  \left(  dx\right)  =\int_{X}\eta^{2}\mu\left(
u,u^{\beta}\right)  \left(  dx\right)  +\int_{X}u^{\beta}\mu\left(  u,\eta
^{2}\right)  \left(  dx\right) \\
\\
=\int_{X}\beta u^{\beta-1}\eta^{2}\mu\left(  u,u\right)  \left(  dx\right)
+\int_{X}u^{\beta}2\eta\mu\left(  u,\eta\right)  \left(  dx\right)  .
\end{array}
\]
Observe that for $\beta\neq-1$, we can write
\[
\frac{4\beta}{\left(  \beta+1\right)  ^{2}}\int_{X}\eta^{2}\mu\left(
u^{\frac{\beta+1}{2}},u^{\frac{\beta+1}{2}}\right)  \left(  dx\right)
=\int_{X}\beta u^{\beta-1}\eta^{2}\mu\left(  u,u\right)  \left(  dx\right)
\]

and
\[
2\int_{X}u^{\frac{\beta-1}2}u^{\frac{\beta+1}2}\eta\mu\left(  \eta,u\right)
\left(  dx\right)  =\frac4{\beta+1}\int_{X}\mu\left(  u^{\frac{\beta+1}2}%
,\eta\right)  u^{\frac{\beta+1}2}\eta\left(  dx\right)  .
\]

Then,
\[
-\int_{X}\beta u^{\beta-1}\eta^{2}\mu\left(  u,u\right)  \left(  dx\right)
\leq2\int_{X}u^{\beta}\eta\mu\left(  \eta,u\right)  \left(  dx\right)  ,
\]
but
\[%
\begin{array}
[c]{c}%
\int_{X}\beta u^{\beta-1}\eta^{2}\mu\left(  u,u\right)  \left(  dx\right)
=\frac{4\beta}{\left(  \beta+1\right)  ^{2}}\int_{X}\eta^{2}\mu\left(
u^{\frac{\beta+1}{2}},u^{\frac{\beta+1}{2}}\right)  \left(  dx\right) \\
\\
\leq\frac{4}{\left(  \beta+1\right)  }\int_{X}\mu\left(  u^{\frac{\beta+1}{2}%
},\eta\right)  u^{\frac{\beta+1}{2}}\eta\left(  dx\right)  .
\end{array}
\]
Taking absolute values gives
\begin{equation}
\frac{\beta}{\left(  \beta+1\right)  }\int_{X}\eta^{2}\left|  \mu\left(
u^{\frac{\beta+1}{2}},u^{\frac{\beta+1}{2}}\right)  \right|  \left(
dx\right)  \leq\int_{X}u^{\frac{\beta+1}{2}}\eta\mu\left(  u^{\frac{\beta
+1}{2}},\eta\right)  \left(  dx\right)  ; \label{b1}%
\end{equation}
recalling the fundamental inequality $2\left|  fg\right|  \left|  \mu\left(
u,v\right)  \right|  \leq f^{2}\mu\left(  u,u\right)  +g^{2}\mu\left(
v,v\right)  ,$ we have
\[
2\left|  fg\right|  \left|  \mu\left(  u,v\right)  \right|  =\left|
u^{\frac{\beta+1}{2}}\eta\right|  \left|  \mu\left(  u^{\frac{\beta+1}{2}%
},\eta\right)  \right|  \leq\frac{\left|  \beta\right|  }{2\left|
\beta+1\right|  }\mu\left(  u^{\frac{\beta+1}{2}},u^{\frac{\beta+1}{2}%
}\right)  \eta^{2}+\frac{\left|  \beta+1\right|  }{2\left|  \beta\right|
}u^{\beta+1}\mu\left(  \eta,\eta\right)  .
\]
Then, from $\left(  \ref{b1}\right)  $ it follows that
\[
\frac{\left|  \beta\right|  }{2\left|  \beta+1\right|  }\int_{X}\eta
^{2}\left|  \mu\left(  u^{\frac{\beta+1}{2}},u^{\frac{\beta+1}{2}}\right)
\right|  \left(  dx\right)  \leq\frac{\left|  \beta+1\right|  }{2\left|
\beta\right|  }\int_{X}u^{\beta+1}\mu\left(  \eta,\eta\right)  \left(
dx\right)  ,
\]
that is
\[
\int_{X}\eta^{2}\left|  \mu\left(  u^{\frac{\beta+1}{2}},u^{\frac{\beta+1}{2}%
}\right)  \right|  \left(  dx\right)  \leq\left(  \frac{\left|  \beta
+1\right|  }{2\left|  \beta\right|  }\right)  ^{2}\int_{X}u^{\beta+1}%
\mu\left(  \eta,\eta\right)  \left(  dx\right)  .
\]
This is the same as $\left(  \ref{a2}\right)  $; beginning from the Sobolev
inequality applied to $u^{\frac{\beta+1}{2}}$ with the same meaning and
definition of the cut-off functions $\eta,$ one gets
\[%
\begin{array}
[c]{c}%
\left(  \frac{1}{m\left(  B_{s}\right)  }\int\limits_{B_{s}}\left|
u^{\frac{\beta+1}{2}}-\bar{u}^{\frac{\beta+1}{2}}{}_{B_{s}}\right|
^{q}m\left(  dx\right)  \right)  ^{\frac{1}{q}}\leq c\tau^{2}c_{1}\left(
sr\right)  sr\left[  \frac{1}{m\left(  B_{s}\right)  }\int\limits_{B_{s}}%
\mu\left(  u^{\frac{\beta+1}{2}},u^{\frac{\beta+1}{2}}\right)  \left(
dx\right)  \right]  ^{\frac{1}{2}}%
\end{array}
\]%
\[%
\begin{array}
[c]{c}%
\leq c\tau^{2}c_{1}\left(  sr\right)  sr\left(  \frac{\left|  \beta+1\right|
}{2\left|  \beta\right|  }\right)  \left(  \int_{B_{s}}u^{\beta+1}\mu\left(
\eta,\eta\right)  \left(  dx\right)  \right)  ^{\frac{1}{2}}%
\end{array}
\]%
\[%
\begin{array}
[c]{c}%
\leq c\tau^{2}\frac{c_{1}\left(  sr\right)  s}{t-s}\left(  \frac{\left|
\beta+1\right|  }{2\left|  \beta\right|  }\right)  \left(  \frac{1}{m\left(
B_{s}\right)  }\int_{B_{s}}u^{\beta+1}m\left(  dx\right)  \right)  ^{\frac
{1}{2}}%
\end{array}
\]%
\begin{equation}
\leq c\tau^{2}\frac{c_{1}\left(  sr\right)  s}{t-s}\left(  \frac{\left|
\beta+1\right|  }{2\left|  \beta\right|  }\right)  \left[  \tau^{2}\left(
\frac{t}{s}\right)  ^{\nu}\right]  ^{\frac{1}{2}}\left(  \frac{1}{m\left(
B_{t}\right)  }\int_{B_{t}}u^{\beta+1}m\left(  dx\right)  \right)  ^{\frac
{1}{2}}. \label{b2}%
\end{equation}
Evaluating the average of $u^{\frac{\beta+1}{2}}$ on $B_{s}$, one gets
\begin{equation}
\underset{B_{s}}{av}u^{\frac{\beta+1}{2}}=\frac{1}{m\left(  B_{s}\right)
}\int_{B_{s}}u^{\frac{\beta+1}{2}}m\left(  dx\right)  \leq\frac{\tau^{2}%
}{m\left(  B_{t}\right)  }\left(  \frac{t}{s}\right)  ^{\nu}\int_{B_{t}%
}u^{\frac{\beta+1}{2}}m\left(  dx\right)  . \label{b3}%
\end{equation}
By H\"{o}lder inequality
\[
\frac{\tau^{2}}{m\left(  B_{t}\right)  }\left(  \frac{t}{s}\right)  ^{\nu}%
\int_{B_{t}}u^{\frac{\beta+1}{2}}m\left(  dx\right)  \leq\left[  \tau
^{2}\left(  \frac{t}{s}\right)  ^{\nu}\right]  ^{\frac{1}{2}}\left[  \frac
{1}{m\left(  B_{t}\right)  }\int_{B_{t}}u^{\beta+1}m\left(  dx\right)
\right]  ^{\frac{1}{2}}%
\]
Putting together $\left(  \ref{b2}\right)  $ and $\left(  \ref{b3}\right)  $,
we see that
\[%
\begin{array}
[c]{c}%
\left(  \frac{1}{m\left(  B_{s}\right)  }\int\limits_{B_{s}}\left|
u^{\frac{\beta+1}{2}}\right|  ^{q}m\left(  dx\right)  \right)  ^{\frac{1}{q}%
}=\left(  \frac{1}{m\left(  B_{s}\right)  }\int\limits_{B_{s}}\left|
u^{\frac{\beta+1}{2}}-\bar{u}^{\frac{\beta+1}{2}}{}_{B_{s}}+\bar{u}%
^{\frac{\beta+1}{2}}{}_{B_{s}}\right|  ^{q}m\left(  dx\right)  \right)
^{\frac{1}{q}}%
\end{array}
\]%
\[%
\begin{array}
[c]{c}%
\leq c\left[  \left(  \frac{1}{m\left(  B_{s}\right)  }\int\limits_{B_{s}%
}\left|  u^{\frac{\beta+1}{2}}-\bar{u}^{\frac{\beta+1}{2}}{}_{B_{s}}\right|
^{q}m\left(  dx\right)  \right)  ^{\frac{1}{q}}+\left(  \frac{1}{m\left(
B_{s}\right)  }\int\limits_{B_{s}}\left|  \bar{u}^{\frac{\beta+1}{2}}{}%
_{B_{s}}\right|  ^{q}m\left(  dx\right)  \right)  ^{\frac{1}{q}}\right]
\end{array}
\]%
\[
\leq\left(  c\tau^{3}\frac{c_{1}\left(  sr\right)  s}{t-s}\frac{\left|
\beta+1\right|  }{2\left|  \beta\right|  }+1\right)  \left(  \frac{t}%
{s}\right)  ^{\frac{\nu}{2}}\left(  \frac{1}{m\left(  B_{t}\right)  }%
\int_{B_{t}}u^{\beta+1}m\left(  dx\right)  \right)  ^{\frac{1}{2}}.
\]
Setting $\beta+1=\tilde{r}$ and $q=2\sigma,$ we see that for any $r$ with
$-\infty<\tilde{r}\leq2,\tilde{r}\neq0,1$%
\[%
\begin{array}
[c]{c}%
\left(  \frac{1}{m\left(  B_{s}\right)  }\int\limits_{B_{s}}u^{\left|
\tilde{r}\right|  \sigma}m\left(  dx\right)  \right)  ^{\frac{1}{\left|
\tilde{r}\right|  \sigma}}%
\end{array}
\]%
\[
\leq\left\{  \left(  c\tau^{3}\frac{c_{1}\left(  sr\right)  s}{t-s}%
\frac{\left|  \tilde{r}\right|  }{2\left|  \tilde{r}-1\right|  }+1\right)
\left(  \frac{t}{s}\right)  ^{\frac{\nu}{2}}\right\}  ^{\frac{2}{\left|
\tilde{r}\right|  }}\left(  \frac{1}{m\left(  B_{t}\right)  }\int_{B_{t}%
}u^{\beta+1}m\left(  dx\right)  \right)  ^{\frac{1}{\left|  \tilde{r}\right|
}}.
\]
We use the iteration argument with any fixed $p$ as a starting value of
$\tilde{r}$, with $-\infty<p\leq2,p\neq0,1.$

\begin{description}
\item [ a)]$p<0$ $\tilde{r}=\sigma^{j}p\rightarrow-\infty$%
\[
\prod_{j=0}^{\infty}\left\{  \left(  c\tau^{3}c_{1}\left(  s_{j}r\right)
a_{j}\frac{\sigma^{j}\left|  p\right|  }{2\left|  \sigma^{j}p-1\right|
}+1\right)  \left(  \frac{s_{j+1}}{s_{j}}\right)  ^{\frac{\nu}{2}}\right\}
^{\frac{2}{\sigma^{j}\left|  p\right|  }}%
\]%

\[
\leq\prod_{j=0}^{\infty}\left\{  \left(  c\tau^{3}c_{1}\left(  s_{j}r\right)
a_{j}\sigma^{j}\left|  p\right|  +1\right)  \right\}  ^{\frac{2}{\sigma
^{j}\left|  p\right|  }}%
\]
Then, after the conversion to the log of the previous quantity we have
\[%
\begin{array}
[c]{c}%
\log\prod_{j=0}^{\infty}\left\{  \left(  c\tau^{3}c_{1}\left(  s_{j}r\right)
a_{j}\sigma^{j}\left|  p\right|  +1\right)  \right\}  ^{\frac{2}{\sigma
^{j}\left|  p\right|  }}=\sum_{j=0}^{\infty}\frac{2}{\sigma^{j}\left|
p\right|  }\log\left\{  \left(  c\tau^{3}c_{1}\left(  s_{j}r\right)
a_{j}\sigma^{j}\left|  p\right|  +1\right)  \right\}
\end{array}
\]%
\[%
\begin{array}
[c]{c}%
=\sum_{j=0}^{\infty}\frac{2}{\sigma^{j}\left|  p\right|  }\left[  \log\left(
c\frac{a_{j}\sigma^{j}}{c_{1}\left(  s_{j}r\right)  \tau^{3}\left|  p\right|
+1}\right)  +\log\left(  \tau^{3}c_{1}\left(  s_{j}r\right)  \left|  p\right|
+1\right)  \right]
\end{array}
\]%
\[%
\begin{array}
[c]{c}%
\leq\frac{2\sigma}{\left(  \sigma-1\right)  \left|  p\right|  }\log\left(
\tau^{3}c_{1}\left(  \frac{r}{2}\right)  \left|  p\right|  +1\right)
+\sum_{j=0}^{\infty}\frac{4}{\sigma^{j}\left|  p\right|  }\left[  \log\left(
c\frac{\left(  j+2\right)  \sigma^{j}}{c_{1}\left(  s_{j}r\right)  \tau
^{3}\left|  p\right|  +1}\right)  \right]
\end{array}
\]%
\[
+\sum_{j=0}^{\infty}\frac{2}{\sigma^{j}\left|  p\right|  }\left[  \log\left(
c\frac{\left(  j+2\right)  \sigma^{j}}{\left(  c_{1}\left(  s_{j}r\right)
\tau^{3}\left|  p\right|  +1\right)  \left(  1-\alpha\right)  }\right)
\right]  \leq\frac{2\sigma}{\left(  \sigma-1\right)  \left|  p\right|  }%
\log\left(  \frac{c\tau^{3}c_{1}\left(  \frac{r}{2}\right)  \left|  p\right|
+1}{\left(  1-\alpha\right)  }\right)
\]

\item[ b)] $0<p<2$%
\[
\prod_{j=0}^{\infty}\left\{  \left(  c\tau^{3}c_{1}\left(  \frac{r}{2}\right)
a_{j}\sigma^{j}p+1\right)  \right\}  ^{\frac{2}{\sigma^{j}p}}\leq\left(
\frac{c\tau^{3}c_{1}\left(  \frac{r}{2}\right)  p+1}{\left(  1-\alpha\right)
}\right)  ^{\frac{2\sigma}{\left(  \sigma-1\right)  p}}%
\]
\end{description}

\begin{lemma}
\label{l3} Let the hypothesis of Lemma $\left(  \ref{l2}\right)  $ hold,
except that now $u$ is a non-negative supersolution of $a\left(  u,v\right)
=0$. For $\frac{1}{2}\leq\alpha<1,$ define $k=k\left(  \alpha,u\right)  $ by
\[
\log k=\frac{1}{m\left(  B_{\alpha}\right)  }\int\limits_{B_{\alpha}}\left(
\log u\right)  dm.
\]
Then for $\lambda>0,$%
\[
m\left(  \left\{  x\in B_{\alpha}:\left|  \log\frac{u\left(  x\right)  }%
{k}\right|  >\lambda\right\}  \right)  \leq\frac{c\tau c_{1}\left(  r\right)
}{\left(  1-\alpha\right)  }\frac{1}{\lambda}m\left(  B_{\alpha}\right)  .
\]
\end{lemma}%

\proof
Letting $\phi=\frac{\eta^{2}}{u},$ it follows that $\left\|  \phi\right\|
_{0}$ is bounded. We observe that
\[%
\begin{array}
[c]{c}%
0\leq\int_{X}\mu\left(  u,\phi\right)  \left(  dx\right)  =\int_{X}\mu\left(
u,\frac{\eta^{2}}{u}\right)  \left(  dx\right)  =
\end{array}
\]%
\[%
\begin{array}
[c]{c}%
\int_{X}\eta^{2}\mu\left(  u,\frac{1}{u}\right)  \left(  dx\right)  +\int
_{X}\frac{1}{u}\mu\left(  u,\eta^{2}\right)  \left(  dx\right)  =\int_{X}%
\eta^{2}\mu\left(  u,u\right)  \left(  -\frac{1}{u^{2}}\right)  \left(
dx\right)  +\int_{X}\frac{1}{u}\mu\left(  u,\eta\right)  2\eta\left(
dx\right)
\end{array}
\]
This implies that%
\[
\int_{X}\eta^{2}\mu\left(  u,u\right)  \left(  \frac{1}{u^{2}}\right)  \left(
dx\right)  \leq2\int_{X}\frac{1}{u}\mu\left(  u,\eta\right)  \eta\left(
dx\right)
\]
Taking the absolute values and recognizing the gradient of the log to the left
side, we have
\[
\int_{X}\eta^{2}\mu\left(  \log u,\log u\right)  \left(  dx\right)  \leq
2\int_{X}\left|  \frac{1}{u}\mu\left(  u,\eta\right)  \eta\right|  \left(
dx\right)  .
\]
Recalling the usual fundamental inequality, we obtain that
\[
2\left|  \frac{1}{u}\eta\right|  \left|  \mu\left(  u,\eta\right)  \right|
\leq\frac{1}{2}\left(  \frac{\eta}{u}\right)  ^{2}\mu\left(  u,u\right)
+2\mu\left(  \eta,\eta\right)  ,
\]
then $\frac{1}{2}\int_{X}\eta^{2}\mu\left(  \log u,\log u\right)  \left(
dx\right)  \leq2\int_{X}\mu\left(  \eta,\eta\right)  \left(  dx\right)  .$ If
we choose as a cut-off function $\eta=1$ on $B_{\alpha}$ s.t. $\left|
\mu\left(  \eta,\eta\right)  \right|  \leq\frac{b}{\left(  1-\alpha\right)
^{2}r^{2}},$ where $b$ is a constant, we get
\[
\int\limits_{B_{\alpha}}\left|  \mu\left(  \log u,\log u\right)  \right|
\left(  dx\right)  \leq\frac{b}{\left(  1-\alpha\right)  ^{2}r^{2}}%
\int\limits_{B_{\alpha}}m\left(  dx\right)  =\frac{bm\left(  B_{\alpha
}\right)  }{\left(  1-\alpha\right)  ^{2}r^{2}}.
\]
By the Poincar\'{e} inequality,
\[%
\begin{array}
[c]{c}%
\int\limits_{B_{\alpha}}\left|  \log u-\overline{\log u}\right|  ^{2}m\left(
dx\right)  \leq cc_{1}^{2}\left(  r\right)  r^{2}\int\limits_{B_{k\alpha}}%
\mu\left(  \log u,\log u\right)  \left(  dx\right)
\end{array}
\]%
\[
\leq bcc_{1}^{2}\left(  r\right)  \frac{m\left(  B_{k\alpha}\right)  }{\left(
1-\alpha\right)  ^{2}}.
\]
It follows that
\[
av_{B_{\alpha}}\left(  \log u\right)  =\frac{1}{m\left(  B_{\alpha}\right)
}\int\limits_{B_{\alpha}}\left(  \log u\right)  m\left(  dx\right)
\rightarrow\log k.
\]
Therefore, by including in $c$ every other inessential constant
\[
\int\limits_{B_{\alpha}}\left|  \log u-\log k\right|  ^{2}m\left(  dx\right)
\leq cc_{1}^{2}\left(  r\right)  \frac{m\left(  B_{k\alpha}\right)  }{\left(
1-\alpha\right)  ^{2}}.
\]
By Chebyshev's inequality, for $\lambda>0,$%
\[%
\begin{array}
[c]{c}%
m\left(  \left\{  x\in B_{\alpha}:\left|  \log\frac{u\left(  x\right)  }%
{k}\right|  >\lambda\right\}  \right)  \leq\frac{1}{\lambda}\int
\limits_{B_{\alpha}}\left|  \log\frac{u}{k}\right|  m\left(  dx\right)
\end{array}
\]%
\[%
\begin{array}
[c]{c}%
\leq\frac{1}{\lambda}\left(  \int\limits_{B_{\alpha}}\left|  \log\frac{u}%
{k}\right|  ^{2}m\left(  dx\right)  \right)  ^{\frac{1}{2}}m\left(  B_{\alpha
}\right)  ^{\frac{1}{2}}\leq\frac{1}{\lambda}\left(  cc_{1}^{2}\left(
r\right)  \frac{m\left(  B_{k\alpha}\right)  }{\left(  1-\alpha\right)  ^{2}%
}\right)  ^{\frac{1}{2}}m\left(  B_{\alpha}\right)  ^{\frac{1}{2}}%
\end{array}
\]%
\[
\leq\frac{1}{\lambda}\tau\left(  \frac{k^{\nu}cc_{1}^{2}\left(  r\right)
}{\left(  1-\alpha\right)  ^{2}}m\left(  B_{\alpha}\right)  \right)
^{\frac{1}{2}}m\left(  B_{\alpha}\right)  ^{\frac{1}{2}}\leq\frac{1}{\lambda
}\tau\left(  \frac{cc_{1}^{2}\left(  r\right)  }{\left(  1-\alpha\right)
^{2}}\right)  ^{\frac{1}{2}}m\left(  B_{\alpha}\right)  .
\]

\begin{lemma}
[Moser]\label{m0}\cite{Moser} Let $c_{1}\left(  r\right)  ,c_{0}\left(
x\right)  $ and $f\left(  x\right)  $ be non-negative bounded functions on a
ball B and in particular let $c_{0}\left(  x\right)  $ belong to
$L_{loc}^{\infty}\left(  X_{0}\right)  $ together with $c_{0}^{-1}\left(
x\right)  $; $c_{1}\left(  r\right)  $ be a decreasing function of $r$ .
Assume that there are constants $D,d$ so that

\begin{description}
\item [(a)]$ess\underset{B_{s}}{\sup}\left(  f^{p}\right)  \leq\left(
\frac{\mu}{t-s}\right)  ^{d}\frac{1}{m\left(  B_{t}\right)  }\int
\limits_{B_{t}}f^{p}m\left(  dx\right)  $ for all $s,t,p$ with $0<p<\frac
{1}{\mu}$ and $\frac{1}{2}\leq s<t\leq1$, where $\mu=\tau^{3}c_{1}\left(
\frac{r}{2}\right)  $

\item[(b)] $m\left(  \left\{  x\in B_{\frac{1}{2}}:\log f\left(  x\right)
>\lambda\right\}  \right)  \leq c\frac{\mu^{\prime}}{\lambda}m\left(
B\right)  ,$ $\forall\lambda>0$, where $\mu=\tau c_{1}\left(  \frac{r}%
{2}\right)  .$
\end{description}

Then, there exists a constant $\gamma$ $\equiv\gamma\left(  c,d\right)  $
s.t.,
\begin{equation}
ess\sup\limits_{B\frac{1}{2}}f\leq\exp\gamma\mu^{2}\mu^{\prime}.\label{m00}%
\end{equation}
\end{lemma}%

\proof
Replacing $f$ with $f_{{}}^{\mu}$ and $\lambda$ with $\lambda\mu$, we simplify
the hypothesis to the case $\mu=1.$ Similarly, we may assume that $m\left(
B\right)  =1$ and the result will be valid for $\mu=1$ too. Define
$\varphi\left(  s\right)  =\underset{B_{s}}{\ \sup}\log f$ for $\frac{1}%
{2}\leq s<1,$ which is a nondecreasing function. Decompose $B_{t}$ into the
sets where $\log f>\frac{1}{2}\varphi\left(  t\right)  $ and where $\log
f\leq\frac{1}{2}\varphi\left(  t\right)  $ and accordingly estimate the
integral
\begin{equation}
\int\limits_{B_{t}}f^{p}m\left(  dx\right)  \leq e^{p\varphi}\frac{2c}%
{\varphi}+e^{p\frac{\varphi}{2}} \label{m1}%
\end{equation}
where $\varphi=\varphi\left(  r\right)  $ and we have used
\begin{equation}
m\left(  \left\{  x\in B_{\frac{1}{2}}:\log f\left(  x\right)  >\lambda
\right\}  \right)  \leq\frac{c}{\lambda} \label{m2}%
\end{equation}
and the normalization $m\left(  B\right)  =1.$ We choose $p$ so that the two
terms on the r.h.s. are equal, i.e.
\begin{equation}
p=\frac{2}{\varphi}\log\left(  \frac{\varphi}{2c}\right)  , \label{m3}%
\end{equation}
provided that this quantity is less than $1$ so that $0<p<\mu^{-1}=1$ holds.
$\frac{2}{\varphi}\log\left(  \frac{\varphi}{2c}\right)  <1$ means that
$c>\frac{\varphi}{2}e^{-\frac{\varphi}{2}},$ but $\frac{\varphi}{2}%
e^{-\frac{\varphi}{2}}$ assumes its maximum at the value $\varphi=2$ and this
means that $\max\left(  \frac{\varphi}{2}e^{-\frac{\varphi}{2}}\right)
=e^{-1}.$ Therefore, if $c>e^{-1}$ then $p<1$; otherwise this requires that
$\varphi=\varphi\left(  r\right)  >c_{1}$ depending only by $c.$ In that case,
we have
\begin{equation}
\int\limits_{B_{t}}f^{p}m\left(  dx\right)  \leq2e^{p\frac{\varphi}{2}}
\label{m4}%
\end{equation}
and hence, by Hp a)
\begin{equation}
\varphi\left(  s\right)  <\frac{1}{p}\log\left(  2\left(  \frac{\mu}%
{t-s}\right)  ^{d}e^{p\frac{\varphi}{2}}\right)  =\frac{1}{p}\log\left(
2\left(  \frac{\mu}{t-s}\right)  ^{d}\right)  +\frac{\varphi\left(  t\right)
}{2} \label{m5}%
\end{equation}
and by $\left(  \ref{m3}\right)  $%
\begin{equation}
\varphi\left(  s\right)  <\frac{\varphi\left(  t\right)  }{2}\left\{
\frac{\log\left(  2\left(  \frac{\mu}{t-s}\right)  ^{d}\right)  }{\log\left(
\frac{\varphi}{2c}\right)  }+1\right\}  . \label{m6}%
\end{equation}
If $\varphi\left(  t\right)  $ is so large that the first term in the
parentheses is less than $\frac{1}{2},$ i.e. if
\begin{equation}
\log\left(  2\left(  \frac{\mu}{t-s}\right)  ^{d}\right)  <\frac{1}{2}%
\log\left(  \frac{\varphi}{2c}\right)  \Longrightarrow\left(  2\left(
\frac{\mu}{t-s}\right)  ^{d}\right)  ^{2}<\frac{\varphi}{2c}\Longrightarrow
\varphi\left(  t\right)  >8c\left(  \frac{\mu}{t-s}\right)  ^{2d} \label{m7}%
\end{equation}
then clearly $\varphi\left(  s\right)  <\frac{3}{4}\varphi\left(  t\right)  .$
Anyway, let us distinguish the case when $p>1$ or $\left(  \ref{m7}\right)  $
fail. This means that:

\begin{enumerate}
\item  if $p>1$, but $\left(  \ref{m7}\right)  $ is still valid, then
$c<\frac{\varphi}{2}e^{-\frac{\varphi}{2}};$ therefore
\[
c<4c\left(  \frac{\mu}{t-s}\right)  ^{2d}e^{-\frac{\varphi}{2}}\Longrightarrow
ce^{\frac{\varphi}{2}}<4c\left(  \frac{\mu}{t-s}\right)  ^{2d},\text{ but
}ce^{\frac{\varphi}{2}}>\varphi\Longrightarrow\varphi<4c\left(  \frac{\mu
}{t-s}\right)  ^{2d}%
\]

\item  if $\left(  \ref{m7}\right)  $ is violated $\Longrightarrow
\varphi\left(  t\right)  <8c\left(  \frac{\mu}{t-s}\right)  ^{2d}.$
\end{enumerate}

In any case $\varphi\left(  t\right)  <8c\left(  \frac{\mu}{t-s}\right)
^{2d}.$ Since $\varphi\left(  s\right)  \leq\varphi\left(  t\right)  ,$ we
have in both cases that
\begin{equation}
\varphi\left(  s\right)  \leq\frac{3}{4}\varphi\left(  t\right)  +\frac
{\gamma_{1}}{\left(  t-s\right)  ^{2d}}, \label{m8}%
\end{equation}
where $\gamma_{1}\equiv\gamma_{1}\left(  c,d\right)  $. For every sequence
$\frac{1}{2}\leq s_{0}\leq s_{1}\leq\ldots\ldots\leq s_{k}\leq1,$ we iterate
the inequality $\left(  \ref{m8}\right)  :$

\begin{enumerate}
\item [ Step 1)]%
\[
\varphi\left(  s_{0}\right)  <\frac{3}{4}\varphi\left(  s_{1}\right)
+\frac{\gamma_{1}\mu^{2}}{\left(  s_{1}-s_{0}\right)  ^{2d}}%
\]
\end{enumerate}

\begin{enumerate}
\item
\[
\varphi\left(  s_{1}\right)  <\frac{3}{4}\varphi\left(  s_{2}\right)
+\frac{\gamma_{1}\mu^{2}}{\left(  s_{2}-s_{1}\right)  ^{2d}},\text{ but }%
\frac{4}{3}\varphi\left(  s_{0}\right)  -\frac{4}{3}\frac{\gamma_{1}\mu^{2}%
}{\left(  s_{1}-s_{0}\right)  ^{2d}}<\varphi\left(  s_{1}\right)  \text{ then}%
\]%
\[%
\begin{array}
[c]{c}%
\frac{4}{3}\varphi\left(  s_{0}\right)  -\frac{4}{3}\frac{\gamma_{1}\mu^{2}%
}{\left(  s_{1}-s_{0}\right)  ^{2d}}<\frac{3}{4}\varphi\left(  s_{2}\right)
+\frac{\gamma_{1}\mu^{2}}{\left(  s_{2}-s_{1}\right)  ^{2d}}\\
\\
\Longrightarrow\varphi\left(  s_{0}\right)  <\left(  \frac{3}{4}\right)
^{2}\varphi\left(  s_{2}\right)  +\gamma_{1}\mu^{2}\left(  \frac{1}{\left(
s_{1}-s_{0}\right)  ^{2d}}+\frac{3}{4}\frac{1}{\left(  s_{2}-s_{1}\right)
^{2d}}\right)
\end{array}
\]
\end{enumerate}

\begin{enumerate}
\item
\[%
\begin{array}
[c]{c}%
\varphi\left(  s_{2}\right)  <\frac{3}{4}\varphi\left(  s_{3}\right)
+\frac{\gamma_{1}\mu^{2}}{\left(  s_{3}-s_{2}\right)  ^{2d}}\\
\\
\Longrightarrow\left(  \frac{4}{3}\right)  ^{2}\varphi\left(  s_{0}\right)
-\gamma_{1}\mu^{2}\left(  \frac{4}{3}\right)  ^{2}\left(  \frac{1}{\left(
s_{1}-s_{0}\right)  ^{2d}}+\frac{3}{4}\frac{1}{\left(  s_{2}-s_{1}\right)
^{2d}}\right)  <\varphi\left(  s_{2}\right) \\
\\
\Longrightarrow\varphi\left(  s_{0}\right)  <\left(  \frac{3}{4}\right)
^{3}\varphi\left(  s_{3}\right)  +\gamma_{1}\mu^{2}\left(  \frac{1}{\left(
s_{1}-s_{0}\right)  ^{2d}}+\frac{3}{4}\frac{1}{\left(  s_{2}-s_{1}\right)
^{2d}}+\left(  \frac{3}{4}\right)  ^{2}\frac{1}{\left(  s_{2}-s_{1}\right)
^{2d}}\right)
\end{array}
\]
\end{enumerate}%

\[
\varphi\left(  s_{0}\right)  <\left(  \frac{3}{4}\right)  ^{k}\varphi\left(
s_{k}\right)  +\gamma_{1}\mu^{2}\sum\limits_{j=0}^{k-1}\left(  \frac{3}%
{4}\right)  ^{j}\frac{1}{\left(  s_{j+1}-s_{j}\right)  ^{2d}}.
\]
By monotonicity, we have $\varphi\left(  s_{k}\right)  \leq\varphi\left(
s_{1}\right)  <\infty$ and letting $k\rightarrow\infty,$ we obtain
\[
\varphi\left(  \frac{1}{2}\right)  \leq\gamma_{1}\mu^{2}\sum\limits_{j=0}%
^{\infty}\left(  \frac{3}{4}\right)  ^{j}\frac{1}{\left(  s_{j+1}%
-s_{j}\right)  ^{2d}}.
\]
The r.h.s. will converge if we choose, for example, $s_{j}=1-\frac{1}{2\left(
1+j\right)  }$%
\[%
\begin{array}
[c]{c}%
\Longrightarrow\varphi\left(  \frac{1}{2}\right)  \leq\gamma_{1}\mu^{2}%
\sum\limits_{j=0}^{\infty}\left(  \frac{3}{4}\right)  ^{j}\left(  \frac
{1}{2\left(  j+2\right)  \left(  j+1\right)  }\right)  ^{2d}=\gamma\mu^{2}\\
\\
\Longrightarrow\quad\underset{B_{\frac{1}{2}}}{\sup}f\leq e^{\gamma\mu^{2}}.
\end{array}
\]%
\ {{\bf Proof of Harnack Inequality.}$\left(
\text{for }\delta>0\right)$}%
Let $\left(  \ref{a1}\right)  ,$ $\left(  \ref{a2}\right)  $ and $\left(
\ref{aa1}\right)  $ hold, $u$ be a non-negative solution of $a\left(
u,v\right)  =0,$ $u\in\mathcal{D}_{loc}\left[  \mathcal{O}\right]  ,\forall
v\in\mathcal{D}_{0}\left[  \mathcal{O}\right]  .$ We wish to apply Lemma
$\left(  \ref{m0}\right)  $ to both $u/k$ and $k/u,$ with $k$ defined by
\begin{equation}
\log k=\frac{1}{m\left(  B_{\frac{3}{2}\alpha}\right)  }\int\limits_{\frac
{3}{2}\alpha}\left(  \log u\right)  m\left(  dx\right)  . \label{h0}%
\end{equation}
Assumption (b) of the lemma holds by applying Lemma $\left(  \ref{l3}\right)
$to $B_{2\alpha}\subset\mathcal{O};$ because of the presence of an absolute
value in the same Lemma and since $\log\left(  u/k\right)  =-\log\left(
k/u\right)  $, assumption (b) holds for both $u/k$ and $k/u.$ Assumption (a)
holds by Lemma $\left(  \ref{l1}\right)  $ applied to $u/k$ with
$d=\frac{\sigma}{\sigma-1}.$ We obtain from $\left(  \ref{m0}\right)  $ both
\[
\underset{B_{\frac{1}{2}}}{ess}\sup\left(  u/k\right)  \leq e^{\gamma\mu
^{2}\mu^{\prime}},\qquad\qquad\underset{B_{\frac{1}{2}}}{ess}\sup\left(
k/u\right)  \leq e^{\gamma\mu^{2}\mu^{\prime}}%
\]
and the result follows by taking the product of these estimates, that is
\[
\underset{B_{\frac{1}{2}}}{ess}\sup u\leq e^{\gamma\mu^{2}\mu^{\prime}%
}\underset{B_{\frac{1}{2}}}{ess}\inf u.
\]%
\ {\bf Proof of Corollary.}%
We may assume without loss of generality that $R\leq R_{0}/4,$ with $R_{0}<1$
and $r\leq\frac{R}{4}$ in such a way that $B\left(  x,4R\right)  \subset$
$B\left(  x,R_{0}\right)  \subset\mathcal{O.}$ Let us define
\[
M_{R}=\underset{B_{R}}{\sup}u\qquad m_{R}=\underset{B_{R}}{\inf}u\qquad
M_{r}=\underset{B_{r}}{\sup}u\qquad m_{r}=\underset{B_{r}}{\inf}u.
\]
Then by applying Harnack inequality to the functions $M_{R}-u,$ $u-m_{R}$ in
$B_{r},$ we obtain
\[
M_{R}-u\leq\underset{B_{r}}{\sup}\left(  M_{R}-u\right)  \leq e^{\gamma
\mu\left(  x,r\right)  }\underset{B_{r}}{\inf}\left(  M_{R}-u\right)
=M_{R}-M_{r}%
\]
and
\[
u-m_{R}\leq\underset{B_{r}}{\sup}\left(  u-m_{R}\right)  \leq e^{\gamma
\mu\left(  x,r\right)  }\underset{B_{r}}{\inf}\left(  u-m_{R}\right)
=m_{r}-m_{R}%
\]
Hence by addition,
\[
M_{R}-m_{R}\leq e^{\gamma\mu\left(  x,r\right)  }\left(  M_{R}-M_{r}%
+m_{r}-m_{R}\right)
\]
so that, writing
\[
\omega\left(  r\right)  =\underset{B_{r}}{osc}u=M_{r}-m_{r}%
\]
we have
\begin{equation}
\omega\left(  r\right)  \leq\left(  1-e^{-\gamma\mu\left(  x,r\right)
}\right)  \omega\left(  R\right)  , \label{c1}%
\end{equation}
with $\omega\left(  R\right)  =\omega\left(  4r\right)  .$ The application of
Lemma 6.5 of \cite{Mosco}, gives the following inequality
\begin{equation}
\omega\left(  r\right)  \leq\exp\left(  -c\int\limits_{\frac{r}{4}}%
^{R}e^{-\gamma\mu\left(  x,\rho\right)  }\frac{d\rho}{\rho}\right)
\omega\left(  R\right)  . \label{c2}%
\end{equation}
Let suppose that $\mu\left(  x,\rho\right)  \approx o\left(  \log\log\frac
{1}{\rho}\right)  $, then
\[%
\begin{array}
[c]{c}%
\omega\left(  r\right)  \leq\exp\left(  -c\int\limits_{\frac{r}{4}}%
^{r}e^{-\gamma\mu\left(  x,\rho\right)  }\frac{d\rho}{\rho}-c\int
\limits_{r}^{R}e^{-\gamma\mu\left(  x,\rho\right)  }\frac{d\rho}{\rho}\right)
\omega\left(  R\right)  ,\\
\\
\Longrightarrow\omega\left(  r\right)  \leq\exp\left(  -\int\limits_{\frac
{r}{4}}^{r}\frac{1}{\log\frac{1}{\rho}}\frac{d\rho}{\rho}-\int\limits_{r}%
^{R}\frac{1}{\log\frac{1}{\rho}}\frac{d\rho}{\rho}\right)  \omega\left(
R\right)  ,
\end{array}
\]
that is,
\begin{equation}
\omega\left(  r\right)  \leq c\frac{\left(  \log\frac{1}{R}\right)  }{\left(
\log\frac{1}{r}\right)  }\omega\left(  R\right)  . \label{c3}%
\end{equation}

\end{document}